\documentstyle{elsart} 

\begin{document}

\newcommand{\ie}{i.e.}
\newcommand{\card}{{\rm Card}}
\newcommand{\nat}{\hbox{I\hskip -2pt N}} 
\newcommand{\pnat}{\hbox{{\scriptsize I}\hskip -2pt {\scriptsize N}}}
\newcommand{\alp}{{\cal A}} 
\newcommand{\balp}{{\cal B}}
\newcommand{\ptf}{X_{\zeta }} 
\newcommand{\sub}{\zeta = (\zeta , \alp , \alpha)}  
\newcommand{\der}{{\cal D}} 
\newcommand{\ene}{{\cal N}} 
\newcommand{\ach}{{\cal H}} 

\def\qed{\vbox{\hrule%
\hbox{\vrule height 5pt depth 0pt\kern 5pt\vrule height 5pt depth 0pt}%
\hrule}}

\begin{frontmatter}
\title{A characterization of substitutive sequences using return words}
\author{F. Durand}
\address{Institut de Math\'ematiques de Luminy, UPR 9016 CNRS, 163 Avenue de Luminy, Case 930 13288 Marseille Cedex 9, France}

\begin{abstract} 
We prove that a sequence is primitive substitutive if and only
if the set of its derived sequences is finite; we defined these sequences here.
\end{abstract}

\end{frontmatter}

\section{Introduction}

The purpose of this paper is to characterize primitive substitutive sequences,
\ie\  sequences defined as the image, under a morphism, of a fixed point of a
primitive substitution. No general characterizations of these sequences are
known. In this paper we give one, it uses a notion we introduce here, the {\em
return word\/}. We don't know reference or study about it, except in the area of
symbolic dynamical systems, where it is closely related to {\em induced
transformations\/} (\cite{denker,peter}). But no combinatoric results are
obtained.

For some classes of sequences arising from substitutions characterizations
exist. The first class appearing in the literature is the class of the sequences
generated by $q$-substitutions, also called {\em uniform tag sequences\/}
(\cite{cob}). These are images, under a "letter to letter" morphism, of a fixed
point of a substitution of constant length $q$.  A. Cobham \cite{cob} proved
that a sequence is generated by a $q$-substitution if and only if it is
$q$-automatic, \ie\ if it is generated by a $q$-automaton (see also
\cite{eilen}). When $q $ is a prime power, an algebraic characterization is
given by G. Christol, T. Kamae, M. Mend\`es France and G. Rauzy \cite{chris}.
And S. Fabre \cite{fabre} generalized the theorem of A. Cobham to a small class
of substitutions of non-constant length, by introducing the concept of
$\Theta$-automaton and $\Theta$-substitution.

It turns out that our characterization is similar to the following known  one:
a sequence $X=(X_n)_{n\in \pnat}$ is generated by a $q$-substitution if and only
if the set of the sequences $ (X_{q^kn+a})_{n\in \pnat}$, where $0\leq a \leq
q^k - 1$ and $k\geq 0 $, is finite (\cite{eilen}).

Our main tool of this paper is the {\em return word over u}, where $u$ is a
non-empty prefix of a minimal sequence $X$. It is defined as a word separating
two successive factors of $u$. We introduce a "coding" of $X$ with the return
words over $u$ and obtain an other minimal sequence ${\cal D}_u(X)$ we call {\em derived sequence\/} of $X$. The main result can be formulated as follows: {\em a sequence is primitive substitutive if and only if the set of derived
sequences is finite}.

The paper is divided into six sections. The next section contains the
definitions, the notation and the basic results needed to state the main
theorem. The third section shows that the condition of the characterization is
sufficient. In the following one, we extend a result, useful for the proof of
the main theorem, due to B. Moss\'e \cite{mosse}, concerning the "power of
words", to primitive substitutive sequences and we deduce some return word
properties of these sequences. In the fifth section we prove the main theorem,
considering first the case of the fixed point, and then the case of its image by
a morphism. The last one deals with perspectives and open questions.

\section{Definitions, notations}

\subsection{Words and sequences}

Let us recall some standard notations we can find in \cite{quef}. We call {\em
alphabet\/} a finite set of elements called {\em letters\/}; symbols $\alp $ and
$\balp $ will always denote alphabets. A {\em word\/} on $\alp$ of {\em length\/} $n$ is an element $x = x_0x_1 \cdots x_{n-1}$ of $\alp^n$, $|x| = n$ is the length of $x$ and $\emptyset$ denotes the {\em empty-word\/} of length 0. The set of non-empty words on  $\alp$ is denoted by $\alp^+$, and $\alp^* = \alp^+ \cup
\{ \emptyset \}$. If ${\cal S}$ is a subset of $\alp^+$, we denote ${\cal S}^+$
the set of words which can be written as a concatenation of elements of ${\cal
S}$, and ${\cal S}^* = {\cal S}^+ \cup \{ \emptyset\}$. A {\em sequence\/} on
$\alp$ is an element of $\alp^{\pnat}$. If $X=X_0X_1\cdots$ is a sequence, and,
$l$ and $k$ are two non-negative integers, with $l\geq k$, we write $X_{[k,l]}$
the word $X_k X_{k+1}\cdots X_{l}$; a word $u$ is a {\em factor\/} of $X$ if there
exists $k\leq l$ with $X_{[k,l]}=u$; if $u=X_{[0,l]}$ for some $l$, we say that
$u$ is a {\em prefix\/} of $X$ and we write $u\prec X$. The empty-word is a prefix
of $X$. The set of factors of length $n$ is written $L_n(X)$, and the set of
factors of $X$ by $L(X)$. If $u$ is a factor of $X$, we will call {\em
occurrence} of $u$ in $X$ every integer $i$ such that $X_{[i,i + |u| - 1]}= u$,
\ie\ such that $u$ is a prefix of $X_{[i,+\infty]}$. The sequence $X$ is {\em
ultimately periodic} if there exist a word $u$ and a word $v$ such that
$X=uv^{\omega}$, where $v^{\omega}$ is the infinite concatenation of the word
$v$. It is {\em periodic\/} if $u$ is the empty-word. We will call $\alp (X)$ the
set of letters occurring in $X$.

If $X$ is a word, we use the same terminology with the similar definitions. A
word $u$ is a suffix of the word $X$ if $|u| \leq |X|$ and $u=X_{[|X|-|u|,
|X|-1]}$. The empty-word is a suffix of $X$.

Every map $\phi : \alp \rightarrow \balp^+$ induces by concatenation a map from
$\alp^* $ to $\balp^*$, and a map from $\alp^{\pnat} $ to $\balp^{\pnat}$. All
these maps are written $\phi $ also.

\subsection{Substitutions and substitutive sequences}

In this part we give a definition of what we call {\em substitution}. In the
literature we can find the notions of {\em iterated morphism\/} (\cite{berstel}),
{\em tag system\/} (\cite{cob}) and {\em CDOL system\/} (\cite{rozen}). They are
similars to the notion of substitution.

\begin{defn}  
A {\rm substitution} is a triple $\zeta = (\zeta , \alp ,
\alpha)$, where $\alp $ is an alphabet, $\zeta $ a map from $\alp $ to $\alp^+$
and $\alpha$ an element of $\alp$ such that: 
\begin{enumerate}  
\item
\label{cond1} 
the first letter of $\zeta (\alpha )$ is $\alpha $, 
\item
\label{cond2} 
$| \zeta^n (\alpha)| \rightarrow +\infty ,$ when $n\rightarrow + \infty $. 
\end{enumerate} 
\end{defn}

This definition could be more general but it is not very restrictive. Indeed, as
we study sequences generated by morphisms we need at least one letter satisfying
the condition \ref{cond2}. Among these letters there is always one, $\beta$,
such that $\beta$ and a power of $\zeta$ satisfy the condition \ref{cond1}.
Details can be found in \cite{quef}.

The sequence $\lim_{i\rightarrow +\infty } \zeta^i (\alpha )$ exists  (see
\cite{quef}) and is called {\em fixed point\/} of $\zeta = (\zeta ,\alp ,
\alpha)$, we  will denote it $X_{\zeta }$; it is characterized by $\zeta(\ptf) =
\ptf$ and $\alpha$ is the first letter of $\ptf$.

\begin{defn} 
We say that a substitution $\zeta = (\zeta , \alp , \alpha)$ is {\rm
primitive} if there exists an integer $k$ such that for all elements $\beta $
and $\gamma $ in $\alp $, $\gamma $ is a factor of $\zeta ^k (\beta )$.
\end{defn}

If $\sub$ is primitive then for all $\beta$ belonging to $\alp$ we have $|
\zeta^n (\beta)| \rightarrow +\infty ,$ when $n$ tends to infinity. 

\begin{defn} 
Let $X$ be a sequence on $\alp$. The sequence $X$ is {\rm minimal}
(or {\rm uniformly recurrent}) if for every integer $l$ there exists an integer
$k$ such that each word of $L_l (X)$ occurs in every word of length $k$.
\end{defn}

Equivalently, a sequence $X$ is minimal if for every non-empty factor $u$ of 
$X$, the maximal difference between two successive occurrences of $u$ is
bounded. If $\zeta$ is primitive then the fixed point $\ptf$ is minimal (see
\cite{quef}). But there exist non-primitive substitutions with minimal fixed
points. Let $\zeta = (\zeta, \{ 1,2,3 \},1)$ be the substitution defined by
$\zeta(1) = 123$, $\zeta(2) = 2$ and $\zeta(3) = 13$. It is not primitive but
its fixed point is minimal.

\begin{defn} 
A sequence $Y$ over the alphabet $\alp$ is called {\rm
substitutive} if there exist a substitution $\zeta = (\zeta , {\cal B} ,
\alpha)$ and a map $\varphi : \cal B \rightarrow \alp $ such that $\varphi
(X_{\zeta }) = Y$, and {\rm substitutive primitive} when $\zeta $ is primitive. 
\end{defn}

It is easy to check that a substitutive primitive sequence is minimal.

\subsection{Return words}

Let $X$ be a minimal sequence over the alphabet $\alp$ and $u$ a non-empty
prefix of $X$. We call {\em return word over u\/} every factor $X_{[i,j-1]}$,
where $i$ and $j$ are two successive occurrences of $u$ in $X$, and ${\cal
H}_{X,u}$ will denote the set of return words over u. Later we will see that
${\cal H}_{X,u}$ is a code;  every element of ${\cal H}_{X,u}^*$ has an unique
decomposition in a concatenation of element of ${\cal H}_{X,u}$ (\cite{beal,berstper}). The sequence $X$ can be written in an unique way as a
concatenation 
\[
X = m_0m_1m_2 \cdots 
\]
of sequence of elements of ${\cal
H}_{X,u}$. Let us give to ${\cal H}_{X,u}$ the linear order defined by the rank
of first appearance in this sequence. This defines a one to one and onto map
\[
\Lambda_{X,u} : {\cal H}_{X,u} \rightarrow \{ 1, \cdots ,\card({\cal
H}_{X,u})\} = {\cal N}_{X,u},
\]
 and the sequence 
\[
\der_u(X) =
\Lambda_{X,u}(m_0)\Lambda_{X,u}(m_1)\Lambda_{X,u}(m_2) \cdots ;
\]
this sequence
of alphabet ${\cal N}_{X,u}$ is called a {\em derived sequence of X\/}; it is
easy to check that $\der_u(X)$ is minimal. Now we are able to state the main
result of this paper:

\begin{thm} 
\label{theo} 
A sequence $X$ is substitutive primitive if and only
if the number of its different derived sequences is finite. 
\end{thm}

We will denote the reciprocal map of $\Lambda_{X,u}$ by $\Theta_{X,u} :
\ene_{X,u} \rightarrow \ach_{X,u}$. The minimal sequence $\der_u(X)$ belonging
to $\ene_{X,u}^{\pnat}$ is characterized by  
\[
\Theta_{X,u}(\der_u(X)) = X. 
\]
The next proposition states some elementary facts about return words, we will
use them very often.

\begin{prop} 
\label{retpref} 
Let $u$ be a non-empty prefix of $X$.
\begin{enumerate} 
\item 
Let $w$ be an element of $\alp^+$. Then $w$ belongs to
$\ach_{X,u}$ if and only if $wu$ is an element of $L(X)$, $u$ is a prefix of
$wu$ and there are exactly two occurrences of $u$ in $wu$. 
\item 
Let
$v_1,v_2,\cdots ,v_k $ be elements of $\ach_{X,u}$. The occurrences of the word
$u$ in $v_1 v_2 \cdots v_k u$ are $0,| v_1|,| v_1|+| v_2| , \cdots ,$
$\sum_{i=1}^k | v_i |$. 
\item 
The set $\ach_{X,u}$ is a code and the map
$\Theta_{X,u} : \ene_{X,u}^*  \rightarrow \ach_{X,u}^*$ is one to one and onto. 
\item 
If $u$ and $v$ are two prefixes of $X$ such that $u$ is a prefix of $v$
then $\ach_{X,v}$ is included in $\ach_{X,u}^+$. 
\item 
Let $v$ be a non-empty
prefix of ${\cal D}_u (X)$ and $w = \Theta_{X,u}(v)u$. Then $w$ is a  prefix $X$
and ${\cal D}_v ({\cal D}_u (X)) = \der_w(X)$.  
\end{enumerate} 
\end{prop} 

\begin{pf} 
Claim 1. comes from the definition of $\ach_{X,u}$. An induction on $k$ prove
the claim 2., and 3. follow from 2..

To prove the claim 4., let $w$ be an element of $\ach_{X,v}$ and $v=ux$. The
word $wv$ is a  factor of $X$, hence the word $wu$ belongs to $L(X)$. The word
$v$ is a prefix of $wv$, consequently the word $u$ is a prefix of $wv$ and the
word $w$ belongs to $\ach_{X,u}^+$. It remains to prove the claim 5..

The word $w$ is a prefix of $X$ because $v$ is a prefix of $\der_u(X)$. Let
$Y=\der_u(X) = q_0q_1q_2\cdots$ and $X= p_0p_1p_2\cdots$ where $q_i$ belongs to 
$\ach_{Y,v}$, $p_i$ belongs to $\ach_{X,w}$ for all integers $i$. Let $i$ be an
element of $\nat$. The word $p_iw$ is a factor of $X$. Moreover $u$ is a prefix
of $w$ and $w$ is a prefix of $p_i w$ then $p_i$ belongs to $\ach^+_{X,u}$. We
write $p_i = \Theta_{X,u}(j_0) \cdots \Theta_{X,u}(j_l)$ where $j_k$ belongs to
$\ene_{X,u}$, $0\leq k\leq l$. 

In addition $v$ is a prefix of $j_0\cdots j_l v$ because $\Theta_{X,u}(v)$ is a
prefix of the word $p_i \Theta_{X,u}(v)$. And the word $j_0\cdots j_l$ belongs
to $\ach_{Y,v}^+$. Let us suppose there are at least three occurrences of $v$ in
$j_0\cdots j_l v$, then there are at least three occurrences of $w$ in $p_i
\Theta_{X,u}(v) u$. It is impossible because $p_i$ is an element of
$\ach_{X,w}$. Then there are exactly two occurrences of $v$ in $j_0\cdots j_l
v$; consequently $\Lambda_{X,u}(p_i) = j_0\cdots j_l $ belongs to $\ach_{Y,v}$.
The sequence $Y$ has an unique concatenation decomposition in elements of
$\ach_{Y,v}$. Hence $\Lambda_{X,u}(p_i) = q_i$, and the map $\Lambda_{X,u}$ from
$\ach_{X,u}$ to $\ene_{X,u}$ is one to one. Then $\Lambda_{Y,v}(q_l) =
\Lambda_{X,w}(p_l)$ for all integers $l$, and finally ${\cal D}_v ({\cal D}_u
(X)) = \der_w(X)$.\qed
\end{pf}

The set $\ach_{X,u}$ is a code, this allows us to extend by concatenation the map
$\Lambda_{X,u}$ to $\ach_{X,u}^*$.

\begin{lem} 
\label{suiteper} 
If the minimal sequence $X$ is ultimately
periodic, then it is periodic. 
\end{lem}

\begin{pf}
There exist two non-empty
words, $u$ and $v$, such that $X = u v^{\omega}$. By minimality, for every
integer $k\geq 1$ there exists an integer $l$ such that the word
$X_{[0,k|v|-1]}$ is a factor of $v^l$. Hence, there exist a suffix, $v_1$, and a
prefix, $v_2$, of $v$ such that the set $\{ k\in \nat ; X_{[0,k |v| - 1]} = v_1
v^{k-1} v_2 \}$ is infinite. We conclude $X = m^{\omega}$ where $m=v_1 v_2$.\qed
\end{pf}

In this paper we only study sequences which are minimals. According to the Lemma \ref{suiteper} we will not make any difference between "periodic" and
"ultimately periodic". We will use "periodic" for both notions.
\begin{prop} 
If the sequence $X$ is periodic and minimal, then there exists a
prefix $u$ of $X$ satisfying: for all words $v$ such that $u$ is a prefix of
$v$ and $v$ is a prefix of $Y$, we have $\card (\ach_{Y,u} )= 1$ and $\ach_{Y,u} = \ach_{Y,v}$. 
\end{prop}

\begin{pf} 
Let $u$ be a word such that $Y=u^{\omega}$.
The word $u$ belongs to $\ach_{Y,u}^+$ and $\Theta_{Y,u}(1)$ is a prefix of $Y$,
consequently $u = \Theta_{Y,u}(1)x$. The word $uu$ is a prefix of $Y$, hence
$\Theta_{Y,u}(1) x \Theta_{Y,u}(1)$ is a prefix of $Y$. Moreover
$\Theta_{Y,u}(1)u$ is a prefix of $Y$, thus $u = x\Theta_{Y,u}(1)$. It follows:
\begin{center}
$Y = u^{\omega} = \Theta_{Y,u}(1) x (\Theta_{Y,u}(1) x )^{\omega} =
\Theta_{Y,u}(1) u^{\omega}$ and $Y = (\Theta_{Y,u}(1))^{\omega}$.
\end{center}

The decomposition of $Y$ in elements of $\ach_{Y,u}$ is unique, consequently
$\ach_{Y,u} = \{ \Theta_{Y,u}(1) \}$.

We know that the word $\Theta_{Y,v}(1)$ belongs to $\ach_{Y,u}^+$. There exists
an integer $k$ such that $\Theta_{Y,v}(1) = (\Theta_{Y,u}(1))^k$. Then $Y =
(\Theta_{Y,v}(1))^{\omega}$ and $\ach_{Y,v} = \{ (\Theta_{Y,u}(1))^k \}$. Hence
the word $\Theta_{Y,u}(1)$ belongs to $\ach_{Y,v}^{+}$. Finally $k=1$, this
completes the proof. \qed
\end{pf}

\section{The condition is sufficient}

To prove that "if the number of the different derived sequences of a minimal
sequence $X$ is finite then $X$ is substitutive primitive". We need the
following proposition. It states that the image of a minimal fixed point by a
morphism from $\alp$ to $\balp^+$ is a substitutive primitive sequence.

\begin{prop}  
\label{imsub}  
Let $\sub $ be a primitive substitution, $\cal B$
an alphabet and $\phi $ a map from $\alp $ to ${\cal B}^+ $. Then the sequence
$\phi (X_{\zeta })$ is substitutive primitive.  
\end{prop} 

\begin{pf}
Let ${\cal C}
= \{(a,k) ; a\in {\alp } \: and \: 1\leq k \leq | \phi (a)| \}$ and $\psi : \alp
\rightarrow {\cal C}^+ $ be the map  defined by:  
\[
\psi (a) = (a,1) \ldots
(a,| \phi (a)| ).
\]
As $\zeta $ is primitive, substituting $\zeta^n$ for
$\zeta$ if needed, we can assume $| \zeta (a)|  \geq| \phi (a)| $, because this
does not change the fixed point.

Let $\tau $ be the map from $\cal C $ to ${\cal C}^+ $ defined by:
\begin{center} 
\begin{tabular}{lllll}     
    & $\tau((a,k))$ & = & $\psi (\zeta (a)_{[k,k]} )$     & if $k < | \phi(a)|
$,\\  and & $\tau((a,| \phi (a)|))$ & = & $\psi (\zeta (a)_{[| \phi (a)|, |
\zeta(a)|]} ) $ & otherwise.\\  
\end{tabular}\\ 
\end{center} 
For $a$ in $\alp $
\begin{eqnarray}
\tau (\psi (a))
 = & \tau ((a,1)\cdots (a,| \phi (a) | ))
 =  \psi (\zeta (a)_{[1,1]} )\cdots \psi (\zeta (a)_{[| \zeta (a) |,| \zeta (a) |]}) \nonumber \\
= & \psi (\zeta (a)) , \nonumber
\end{eqnarray}
thus  $ \tau (\psi (X_{\zeta}))=\psi (\zeta
(X_{\zeta}))=\psi(X_{\zeta }).$

In this way $\psi (X_{\zeta })$ is the fixed
point of $\tau = (\tau , {\cal C}, (\alpha ,1))$ and $\psi (X_{\zeta })=X_{\tau
}$. The substitution $\zeta $ is primitive and $\tau^n \psi =\psi \zeta^n$,
hence $\tau$ is primitive. If $\chi $ is the map from $\cal C$ to $\cal B$ which
sends the $k$-th letter of $\phi(a)$ to $(a,k)$, we obtain~:  
\[
\chi (\psi (a))
= \chi ((a,1) \cdots (a,| \phi (a) | ))= \phi (a),
\]
and $\chi (X_{\tau}) =
\chi (\psi (X_{\zeta})) = \phi (X_{\zeta }).$

Therefore $\phi (X_{\zeta })$ is
primitive substitutive.\qed
\end{pf}

\begin{lem} 
\label{mrcroi} 
Let $X$ be a minimal sequence, on the alphabet
$\alp$, which is not periodic. Then 
\[
m_n = \inf \{ | v | ; v\in \ach
_{X,X[0,n]} \}  \rightarrow +\infty \:\: when \:\: n\rightarrow +\infty. 
\]
\end{lem} 

\begin{pf}
We have seen, it is the Proposition \ref{retpref}, that 
$\ach_{X,X[0,n+1]}$ is included in $\ach_{X,X[0,n]}^+$, then $m_{n} \leq
m_{n+1}$. Let suppose $(m_n)_{n\in \pnat}$ stationary at the rank $n_0$: there
exist an integer $k$ and, for every $n > n_0 $, a word $v_n$ where $| v_n| = k$
and $v_n$ is an element of $\ach_{X,X[0,n]}$. If $n\geq k$ then $X[0,n]$ is a
prefix of $v_n X[0,n]$. Therefore, for all integers $j$ such that $0\leq j \leq
n-k$, we deduce $X[j] = X[k+j]$. It follows that $X$ is periodic with period
$k$. This completes the proof.\qed
\end{pf}

In the following proposition we choose two prefixes $u$ and $v$ of a
non-periodic minimal sequence $X$ satisfying: $u$ is a prefix of $v$, such that
each word $tu$, where $t$ belongs to $\ach_{X,u}$, is a factor of every $w$
belonging to $\ach_{X,v}$. The minimality of $X$ and the Lemma \ref{mrcroi}
allows us to set such a hypothesis.

\begin{prop}  
\label{rproj}  
Let $X$  be a non-periodic minimal sequence. Let
$u$ and $v$ be two prefixes of $X$, where $u$ is a prefix of $v$, such that each
word $tu$, where $t$ belongs to $\ach_{X,u}$, is a factor of every $w$ belonging
to $\ach_{X,v}$. If $\der_u (X) =\der_v (X)$, then $\der_u(X)$ is the fixed
point of a primitive substitution and $X$ is substitutive primitive.  
\end{prop}   

\begin{pf}
Suppose $\der_u (X) =\der_v (X)$. Then we have 
\[
\ene_{X,u} = \alp (\der_u (X)) = \alp (\der_v (X)) = \ene_{X,v}.
\]
And ${\cal H}_{X,v}$ is included in $ {\cal H}_{X,u}^+$ because $u$ is a prefix
of $v$. Therefore we obtain $\card ( \ene_{X,u} ) =  | \{ \Lambda_{X,u}
\Theta_{X,v} (a); a\in \ene_{X,v} \} | $. Let us consider the map $\zeta$ from
$\ene_{X,u}$ to $\ene_{X,u}^+$ defined by $\zeta = \Lambda_{X,u} \Theta_{X,v} $.
It is easy to see that the first letter of $\zeta(1)$ is 1. Moreover, for every
$i,j$ belonging to $\ene_{X,u}$, $\Theta_{X,u}(j)u$ is a factor of
$\Theta_{X,v}(i)$. Therefore $j$ appears in $\Lambda_{X,u} \Theta_{X,v}(i)=\zeta
(i)$ and $\zeta = (\zeta , \ene_{X,u}, 1)$ is a primitive substitution. And
$\der_u (X)$ is the fixed point of $\zeta$ because 
\begin{center}
$
\zeta (\der_u(X)) =
\Lambda_{X,u} \Theta_{X,v} (\der_u(X)) = \Lambda_{X,u} \Theta_{X,v} (\der_v(X))
= \Lambda_{X,u}(X) = \der_u (X).
$
\end{center}
With the Proposition \ref{imsub}, we conclude that $X$ is substitutive primitive because $\Theta_{X,u}(\der_u(X))$ $=X$.\qed
\end{pf}

Now we are able to show that the condition of the characterization is sufficient.
\begin{thm} 
If the minimal sequence $X$ has a finite number of derived
sequences, then  it is substitutive primitive. 
\end{thm}

\begin{pf}
If $X$ is
periodic, then it is easy to check that $X$ is the fixed point of a primitive
substitution of constant length.

Let $X$ be non-periodic. There exists a prefix $u$ of $X$ such that the set $K =
\{ v\prec X; \der_v(X) = \der_u(X) \}$ is infinite. By minimality, we can choose
$n$ so large that every factor of length $n$ of $X$ contains each element of
$\ach_{X,u}$ as a factor. By Lemma \ref{mrcroi} there exists a word $v$
belonging to $K$ such that $|w|\geq n$ for all $w$ belonging to $\ach_{X,v}$.
With $u$ and $v$ the hypothesis of Proposition \ref{rproj} are fulfilled, this
completes the proof.\qed
\end{pf}

\section{Return words of a substitutive primitive sequence}

In this section $\sub$ is a primitive substitution, $\phi : \alp \rightarrow
{\cal B}$ is a map, we call it {\em projection\/}, and $Y=\phi (\ptf)$. We have to
introduce some notations. We define 
\[
S(\zeta) = \sup \{ | \zeta(a) | ; a\in
\alp \} \:\: {\rm and} \:\: I(\zeta) = \inf \{ | \zeta(a) | ; a\in \alp \}.
\]

\subsection{Power of words in a substitutive primitive sequence}

The following result comes from \cite{quef}. 
\begin{lem} 
\label{submaj} 
There exists a constant $Q$ such that for all integers $n$ 
\[
S(\zeta^n) \leq Q
I(\zeta^n).
\]
\end{lem}

\begin{defn} 
A word $v$ is {\rm primitive} if it does not exist an integer $n$
and a word $w\not = v$ such that $v = w^n$.  
\end{defn}

A word is always a power of a primitive word. The proof of the next result can
be found in \cite{mosse}, so we omit the proof.

\begin{prop} 
\label{power} 
Let $w = v^n$, where v is a primitive word and $n\geq
2$. If $vuv$ is a factor of $w$, then $u$ is a power of $v$. 
\end{prop}

The two following results have been proved by B. Moss\'e in \cite{mosse} for
fixed points of primitive substitutions, we adapted her proof to substitutive
primitive sequences.

\begin{lem} 
\label{subper} 
If there exist a primitive word $v$ and two integers
$N$ and $p$ such that:
\begin{enumerate} 
\item 
for all words $ab$ of $L_2(\ptf
)$, $\phi(\zeta^p(ab))$ is a factor of $v^N$, 
\item 
$2| v| \leq I(\zeta^p)$,
\end{enumerate} 
then $Y$ is periodic. 
\end{lem}

\begin{pf}
According to the conditions 1. and 2., for all letters $a$ of $\alp$,
there exist an integer $n_a\geq 1$, a prefix $w_a$ of $v$ and a suffix $v_a$ de
$v$ such that $|w_a|,|v_a| < |v|$ and $\phi (\zeta^p (a)) = v_a v^{n_a}w_a$. If
$ab$ belongs to $L_2(\ptf)$ then $v^{n_a}w_av_b v^{n_b}$ belongs to $L(Y)$. The
word $v$ is primitive then $w_av_b=v$ or $w_av_b=\emptyset$. Hence
$Y=\phi(\ptf)=\phi(\zeta^p (\ptf))$ is periodic.\qed
\end{pf}

\begin{thm} 
\label{puismot} 
If $Y$ is not periodic, there exists an integer $N$
such that $w^N$ is a factor of $Y$ if and only if $w=\emptyset$. 
\end{thm}

\begin{pf}
By minimality, there exists an integer $r$ such that each word of length
2 of $Y$ is a factor of each factor of length $r$ of $Y$. Let $v$ be a non-empty
primitive word of $\balp^+ $, $M>0$ be an integer such that $v^M$ belongs to
$L(Y)$ and $p$ be the integer defined by $I(\zeta^{p-1}) \leq 2| v|
<I(\zeta^p)$. There is an occurrence of each word $\phi (\zeta^p (ab))$, where
$ab$ belongs to $L(\ptf)$, in each word of length $2\sup \{ | \zeta^p(w)| ; w\in
L(\ptf), | w| = r \}$. Hence, according to the Lemma \ref{subper}, we deduce  
\[
| v^M | = M| v| < 2 \sup \{ | \phi(\zeta^p (w)) | ; w\in L(\ptf), | w| =
r\}.
\]
Then, using the Lemma \ref{submaj}, and its constant $Q$, 
\begin{eqnarray}
M 
< & \frac{2
\sup \{ | \phi(\zeta^p (w)) | ; w\in L(\ptf), | w| =
r\}}{\frac{1}{2}I(\zeta^{p-1}) } \nonumber \\
\leq & 4r\frac{S(\zeta^p)}{I(\zeta^{p-1} )} 
\leq 4rS(\zeta)\frac{S(\zeta^{p-1})}{I(\zeta^{p-1})} 
\leq 4rS(\zeta)Q. \nonumber
\end{eqnarray}
With $N = 4rS(\zeta)Q$ we complete the proof.\qed
\end{pf}

\subsection{Return words}

Here we prove that, for a substitutive primitive sequence, the length of the
return words over $u$ is proportional to the length of $u$, and that the number
of return words are bounded independently of $u$.

\begin{thm} 
\label{encadr}
If $Y$ is not periodic, there exist three positive
constants $K$, $L$ and $M$ such that: for all non-empty prefixes $u$ of $Y$,
\begin{enumerate} 
\item 
for all words $v$ of $\ach_{Y,u}$, $L| u| \leq | v| \leq
M| u|$, 
\item 
$\card (\ach_{Y,u})\leq K$. 
\end{enumerate} 
\end{thm}

\begin{pf}
Let
$u$ be a non-empty prefix of $Y$.\\ 1. Let $v$ be an element of $\ach_{Y,u} $
and $k$ be the smallest integer such  that $| u| \leq I( \zeta^k)$. The choice
of $k$ entails that there exists an element $ab$ of $L_2(\ptf)$ such that $u$ is
a factor of $\phi (\zeta^k (ab))$. Let $R$ be the largest difference between two
successive occurrences of an element of $L_2(\ptf)$ in $\ptf$. We have  
\[
| v| \leq R S(\zeta^k) \leq RQ I(\zeta^k) \leq RQ S(\zeta) I(\zeta^{k-1}) \leq RQ
S(\zeta) | u| ,
\]
and we put $M=RQ S(\zeta) $.

Theorem \ref{puismot} gives us $| v| \geq | u| / N$, where $N$ is the constant
of this theorem.\\ 2. Let $n$ be the smallest integer such that $I(\zeta^n) \geq
(M+1)| u |$. For all $v$ belonging to $\ach_{Y,u}$ we have $| vu| \leq (M+1) |
u| $. Hence there exists a word $ab$ of $L_2 (\ptf)$ such that $vu$ is a factor
of $\phi (\zeta^n (ab))$. The word $\phi (\zeta^n (ab))$ contains at the most $q
= | \phi (\zeta^n (ab))|  L^{-1} | u|^{-1}$ occurrences of $u$. Consequently at
the most $q$ factors $vu$ where $v$ belongs to $\ach_{Y,u}$. On the other hand
\begin{eqnarray}
| \phi (\zeta^n (ab))| 
\leq & 2S(\zeta^n) 
\leq  2S(\zeta)S(\zeta^{n-1}) \nonumber \\
\leq & 2S(\zeta)QI(\zeta^{n-1}) 
\leq  2S(\zeta)Q(M+1) |u| ,\nonumber
\end{eqnarray}
hence $q\leq
2S(\zeta)Q(M+1)L^{-1}$ and   
\[
\card (\ach_{Y,u})\leq 2S(\zeta)Q(M+1) (\card
(\alp))^2 L^{-1},
\]
this completes the proof.\qed
\end{pf}

\section{The condition is necessary}

To prove "if $X$ is a primitive substitutive sequence then the number of its
different derived sequences is finite", this section is divided into two
subsections. One is devoted to the case of the fixed point of a primitive
substitution, and the other to the general case.

\subsection{The case of the fixed point of a primitive substitution}

\begin{prop}
\label{ptfixe=TR} 
Let $\sub $ be a primitive substitution and $u$ be
a non-empty prefix of $\ptf $. The derived sequence $\der_u (\ptf) $ is the
fixed point of a primitive substitution $\tau_{\ptf,u} = (\tau_{\ptf,u} ,
\ene_{\ptf,u}, 1)$.  
\end{prop} 

\begin{pf} 
Let $i$ be an element of $\ene_{\ptf,u}$
and $v = \Theta_{\ptf,u}(i)$. The word $v$  belongs to $\ach_{\ptf,u}$, hence
$vu$ is a factor of $\ptf$ and $\zeta(v)\zeta(u)$ too. The word $u$ is a prefix
of $\ptf$ and $u$ is a prefix of $vu$. Therefore $u$ is a prefix of $\zeta (u)$
and $\zeta (v) \zeta (u)$. And finally $\zeta (v)$ belongs to $\ach_{\ptf,u}^+$.
So, we can define the map $\tau_{\ptf,u}$ from $\ene_{\ptf,u}$ to
$\ene_{\ptf,u}^+$ by:  
\[
\tau_{\ptf,u} (i) = \Lambda_{\ptf , u} (\zeta
(\Theta_{\ptf,u}(i))),\:\: i\in \ene_{\ptf,u}.
\]
In this way we have  
\begin{eqnarray}
\tau_{\ptf,u} ( \der_u(\ptf)) & 
= \Lambda_{\ptf , u} (\zeta (\Theta_{\ptf,u} (\der_u(\ptf)))) 
= \Lambda_{\ptf , u}(\zeta (\ptf))
=\Lambda_{\ptf , u}(\ptf ) \nonumber \\
& =\der_u(\ptf). \nonumber
\end{eqnarray}
The map $\tau_{\ptf,u} = (\tau_{\ptf,u}, \ene_{\ptf,u},1)$
is a substitution because $1$ is a prefix of $\tau_{\ptf,u} (1)$ and
\[
\lim_{n\rightarrow +\infty} | \tau_{\ptf,u}^n (i)| = \lim_{n\rightarrow
+\infty} | \Lambda_{\ptf ,u}(\zeta^n (\Theta_{\ptf,u}(i)))| = +\infty, \forall
i\in \ene_{\ptf,u}.
\]
And $X_{\tau_{\ptf,u}} = \der_u(\ptf)$ is its fixed
point. Finally, it is primitive because $\zeta$ is primitive.\qed
\end{pf}

\begin{thm}
\label{suitder} 
Let $\sub$ be a primitive substitution. The number
of dif\-ferent deri\-vated sequences of $\ptf$ is finite. 
\end{thm} 

\begin{pf}
Let
$\ptf $ be non-periodic. Let $i$ be an element of $\ene_{\ptf,u}$ and $v =
\Theta_{\ptf,u}(i)$ an element of $\ach_{\ptf ,u}$. We take the notations of
Theorem \ref{encadr}. We have $| v| \leq M| u|$ and $| \zeta (v) | \leq S(\zeta)
M | u|$. The length of each element of $\ach_{\ptf,u}$ is larger than $L| u|$.
Then we can decompose $\zeta(v)$ in at the most $S(\zeta)ML^{-1}$ elements of
$\ach_{\ptf,u}$, so $| \tau_{\ptf,u}(i) | \leq S(\zeta)ML^{-1}$. 

There exists $K$ such that $\card(\ach_{\ptf,u}) \leq K$ for all non empty
prefixes of $\ptf$. Therefore we can deduce there is a finite number of
alphabets $\ene_{\ptf,u}$, and a finite number of substitutions $\tau_{\ptf,u} =
(\tau_{X,u} , \ene_{\ptf,u}, 1)$~: the set $\{ \der_u(\ptf) ; u\prec \ptf \}$ is
finite.

Let $\ptf $ be periodic. There exists a prefix $u$ of $\ptf$ satisfying $\card
(\ach_{\ptf ,v} ) = 1$ for all words $v$ such that $u$ is a prefix of $v$ and
$v$ is a prefix of $\ptf$. Then $\der_v (\ptf) = 1^{\omega}$ for all words $v$
such that $u$ is a prefix of $v$ and $v$ is a prefix of $\ptf$. The proof is
completed.\qed
\end{pf}

\subsection{The case of a substitutive primitive sequence}

Here we end the proof of the main theorem.

\begin{thm} 
Let $Y$ be a substitutive primitive sequence, the number of its
dif\-ferent deri\-vated sequences is finite. 
\end{thm}

\begin{pf}
If $Y$ is
periodic, we treated this case in the proof of Theorem \ref{suitder}.

Let $Y$ be non-periodic. We have $Y=\phi (\ptf)$ where $\phi $ is a projection
and $\zeta $ a primitive substitution. We proved that there exist three positive
constants $K$, $L$ and $M$ such that for all non-empty prefixes $u$ of $Y$, and
for all words $v$ of  $\ach_{Y,u}$ we have $L| u| \leq | v| \leq M| u|$ and
$\card (\ach_{Y,u})\leq K$.

Let $v$ be a non-empty prefix of $Y$, and $u$ the prefix of $\ptf$ of length
$|v|$. We have $v = \phi (u)$. Let $i$ be an element of $\ene_{\ptf ,u}$ and $w
= \Theta_{\ptf ,u}(i)$ an element of $\ach_{\ptf,u}$. As $w$ belongs to
$L(\ptf)$ the word $\phi (w) v $ belongs to $L(Y)$. Moreover $u$ is a prefix of
$wu$, then $v$ is a prefix of $\phi (w) v$. Therefore $\phi (w)$ belongs to
$\ach^{+}_{Y,v}$ and this word has an unique concatenation decomposition in
elements of $\ach_{Y,v}$. We can define the map $\lambda_u : \ene_{\ptf,u}
\rightarrow \ene^*_{Y,v}$ by $\lambda_u = \Lambda_{Y,v}  \phi \Theta_{\ptf,u}$.
So, we have $\lambda_u(\der_u (\ptf)) = \der_v(Y)$, $\Theta_{Y,v}  \lambda_u  =
\phi \Theta_{\ptf ,u}$ and 
\begin{eqnarray}
|\lambda_u (i) | L |v| 
\leq & |\lambda_u (i) |\inf \{ |w|;w\in \ach_{Y,v} \} 
\leq  | \Theta_{Y,v}(\lambda_u (i))| \nonumber \\
 = & | \Theta_{\ptf,u}( i)| 
\leq \sup \{ |w|;w\in \ach_{X,u} \} 
\leq M|u| \leq M|v|.\nonumber
\end{eqnarray}
Hence $|\lambda_u (i) |\leq M/L.$

The sets $\{ \lambda_u ; u\prec \ptf
\}$ and $\{ \der_u(\ptf ); u\prec \ptf \}$ are finites. Therefore the set $\{
\der_v(Y ); v\prec Y \} =  \{  \lambda_u (\der_u(\ptf )); u\prec \ptf, |u| = |v|
\}$ is finite. This completes the proof.\qed
\end{pf}

Theorem \ref{theo} is proved.

\section{Open problems and perspectives}

\subsection{Topological dynamical system}

Let $X$ be a minimal sequence and $\Omega = (\overline{\{ T^n X ; n\in \nat
\}},T)$, where $T$ is the shift, be the dynamical system generated by $X$ (see \cite{denker,peter} or \cite{quef}). In dynamical topological terms, we
can formulate Theorem \ref{theo} as follows: $X$ is primitive substitutive if
and only if the number of induced systems on cylinder generated by a prefix is
finite. We don't know whether it is true if we induce on any cylinders, or on
any clopen sets.

Our method characterizes only one point of $\Omega$, is it possible to characterize all the points belonging to $\Omega$ by analysis methods?

\subsection{Minimality and primitivity} 

Our method allows us to characterize the primitive substitutive sequences, but
we do not know if every substitutive minimal sequence is primitive, \ie\
generated by a primitive substitution. We can only prove this for the minimal
fixed points arising from non-primitive substitutions. The following example is
significant and give a sketch of the proof.

Let $\zeta = (\zeta , \{ 1,2\}, 1)$ be the non-primitive substitution defined by $\zeta(1) = 1211$ et $\zeta (2) = 2$. The sequence $\ptf$ is minimal. There
exists a prefix $u$ of $\ptf$ such that $\der_u(\ptf)$ is the fixed point of the primitive substitution $\tau_{\ptf , u}$. Then $\ptf$ is primitive substitutive
because $\Theta_{\ptf,u}(\der_u(\ptf)) = \ptf$. The sequence $\ptf$ is the image
under $\varphi : \{ 1,2,3\} \rightarrow \{ 1,2\}$, defined by $\varphi
(1)=\varphi (3)=1$ and $\varphi (2)=2$, of the fixed point of the primitive
substitution $\tau = (\tau ,\{ 1,2,3\}, 1)$, defined by $\tau (1) = 12$, $\tau
(2) = 312$ and $\tau (3) = 1233$. The proof of Theorem \ref{imsub} gives us the
morphisms $\varphi$ and $\tau$.

It remains to treat the case of the substitutive minimal sequences arising from
sequences which are not minimals.

\subsection{Complexity, return words and $S$-adic sequences}

B. Moss\'e recently proved that if $X$ is the fixed point of a primitive 
substitution then the sequence $(p_X(n+1) - p_X(n))_{n\in \pnat}$ is bounded,
where $p_X(n)$ is the number of the words of length $n$ belonging to $L(X)$.
With our results we are able to give a short proof of this property, and to
extend it to primitive substitutive sequences. Indeed we can prove, with the
help of the graph of words of Rauzy \cite{rauzy}, the following property: \\ if
$Y$ is a minimal sequence such that $({\rm Card}({\cal H}_{Y,Y_{[0,n]}}))_{n\in
\pnat}$ is bounded then $(p_Y(n+1) - p_Y(n))_{n\in \pnat}$ is bounded. 

A $S$-adic sequence $Y$ on the alphabet $\alp$ is given by an infinite product
of substitutions belonging to a finite set ${\cal S}$. B. Host formulated the
following conjecture: a minimal sequence $Y$ is $S$-adic if and only if
$(p_Y(n+1) - p_Y(n))_{n\in \pnat}$ is bounded. 

If the matrices of the substitutions of ${\cal S}$ (see \cite{quef}) are
strictly positives then we can prove that $({\rm Card}({\cal
H}_{Y,Y_{[0,n]}}))_{n\in \pnat}$ is bounded, by extending the proof of the
primitive substitutive case. It follows that $(p_Y(n+1) - p_Y(n))_{n\in \pnat}$
is bounded.

\end{document}